\documentclass[12pt]{article}

\nonstopmode

\usepackage[left=1.5cm,right=1.5cm,top=2cm,bottom=2cm]{geometry}
\usepackage[utf8]{inputenc}
\usepackage[english]{babel}
\usepackage[T1]{fontenc}
\usepackage{amsmath}
\usepackage{amsthm}
\usepackage{amsfonts}
\usepackage{amssymb}
\usepackage{color}
\usepackage{multicol}
\usepackage{times}
\usepackage{setspace}
\usepackage{hyperref}
\usepackage{hyperref}
\hypersetup{
colorlinks=true, 
breaklinks=true, 
urlcolor= blue, 
linkcolor= blue, 
citecolor=blue, 
}
\usepackage{graphicx}
\usepackage{lscape}
\usepackage{setspace}

\newtheorem{theorem}{Theorem}[section]
\newtheorem{lemma}{Lemma}[section]

\newtheorem{cor}{Corollary}[section]

\newtheorem{definition}{Definition}[section]

\newcommand{\bZ}{\ensuremath{\mathbb{Z}}}
\newcommand{\bQ}{\ensuremath{\mathbb{Q}}}

\newcommand{\bL}{\ensuremath{\mathbb{L}}}

\newcommand{\klammern}[4][]%
{\ifthenelse{\equal{#1}{}}{\left#2}{\csname#1\endcsname#2}%
#4\ifthenelse{\equal{#1}{}}{\right#3}{\csname#1\endcsname#3}}

\title{The Diophantine equation $\displaystyle a \left(\frac{b^k - 1}{b - 1}\right) = \mathcal{U}_n - \mathcal{U}_m$} 
\author{P. Tiebekabe$^{1,2}$, I. Diouf$^{2}$, A. Tall$^{2}$, and K. R. Kakanou$^{2}$}
\date{}
\begin{document}
\maketitle
\begin{abstract}
\noindent Here, we find all positive integer solutions of the Diophantine equation in the title, where $(\mathcal{U}_n)_{n\geqslant 0}$  is the generalized Lucas sequence $\mathcal{U}_0=0, \ \mathcal{U}_1=1$  and $\mathcal{U}_{n+1}=r \mathcal{U}_n +s \mathcal{U}_{n-1}$ with $r$ and $s$ integers such that $\Delta = r^2 +4 s >0$.





\end{abstract}
\textbf{Keywords}:  Linear forms in logarithm; Diophantine equations; Linear recurrent sequences;
		Lucas  number; Pell number, Baker's theory, Reduction procedure\\
\textbf{2020 Mathematics Subject Classification: 11B39, 11J86, 11D61, 11Y50.}\\
\textbf{Affiliations:}\\
$^1$ Universi\'e de Kara, Facult\'e des Sciences et Techniques (FST), D\'epartement de Math\'ematiques, Kara, Togo.\\
$^{2}$ Universit\'e Cheikh Anta Diop de Dakar (UCAD), Faculté des Sciences et Techniques (FST), Laboratoire d’Algèbre, de Cryptologie, de G\'eom\'etrie Alg\'ebrique et Applications (LACGAA), Dakar, BP.  Fann 5005 Dakar-Senegal.

\section{Introduction}\label{sec1}
Numerical patterns and relationships have intrigued mathematicians for centuries. From prime numbers to Fibonacci sequences, these patterns often hold a captivating allure. In this article, in one hand, we delve into an intriguing connection between b-repdigits and generalized Lucas numbers, unearthing an elegant representation of b-repdigits as the difference between two generalized Lucas numbers. A $b$-repdigit refers to a number composed of a repeated digit in base $b$. On the other hand, generalized Lucas numbers, denoted as $\mathcal{U}_n$, are a sequence of numbers that exhibit a recursive pattern, making them a fascinating topic of study in their own right. Recent papers have contributed significantly to understanding repdigits, exploring various aspects of these intriguing numerical patterns. Studies have focused on their properties, connections with recurrent sequences, and representations as concatenations, differences, or products of other mathematical entities. For details of recent work related to this problem studied in this paper, we refer the reader to \cite{erd1, adedji, adedji2, siar, erd2, bhoi}. They contribute to the field of number theory, inspiring further research in the exploration of repdigits and their intricate connections with other mathematical entities. With this in mind, we have studied b-repdigits, which are the difference between two generalized Lucas numbers, and given an application to the case of the Pell sequence in the decimal base.
	
 \noindent Our paper is organized as follows, in the first section we give our main results. The second one is for recalls, the third for the proof of the main results, and the last to present the application of the fundamental theorem is the particular case of Pell numbers in decimal base.
	
	\section{Results}\label{sec2}
	
	The main results of this paper are the following.
	\begin{theorem}\label{main}$\\$
		Let	$b$ a positive integer such that $b\geq 2$. If $k, m$ and $n$ are positive integers that satisfy the Diophantine equation $a\left( \dfrac{b^k-1}{b-1}\right)=\mathcal{U}_n-\mathcal{U}_m$ with $n>m$ and $1\leqslant a\leqslant b-1$, then 
		$$
		k<1+ n\dfrac{\log \delta}{\log b}
		$$
		and 
		$$
		n\log \delta -\log(8.1\sqrt{\Delta})< 2\cdot 10^{12}(1+\log D)\log \delta\log b\cdot \xi$$   where $$\xi= \log(4 b^2 \Delta (1+3 \sqrt{\Delta}) +2\cdot 10^{12} (1+\log D)  \cdot \log \delta \cdot  \log b \cdot (2 \log b + \log \Delta) $$ with $$\quad D=1+ n\dfrac{\log \delta}{\log b}.$$
	\end{theorem}
	Moreover, the above result implies.
	\begin{cor}$\\$
	The Diophantine equation (\ref{Eq1}) has only finitely many solutions in positive integer $k, m, n, b$ and $a$.
	\end{cor}
	By considering the cases $b=10$ and the particular Pell numbers, we get the following result.
	\begin{theorem}\label{theo}$\\$
		The only repdigits which are differences between two Pell numbers are :
		\[1, 3, 4, 7, 11\  \text{and} \  99.\]
		Moreover, we have \\ 
		\begin{center}
			\begin{tabular}{|c|c|c|}
				\hline
				n& m & $P_n-P_m$\\
				\hline
				2 & 1 & 1\\
				\hline
				3 & 2 & 3 \\
				\hline
				3 & 1 & 4 \\
				\hline
				4 & 3 & 7 \\
				\hline
				4 & 1 & 11 \\
				\hline
				7 & 6 & 99 \\
				\hline
			\end{tabular}
		\end{center}

		where $P_l$ is $l$-th term of Pell sequence.
	\end{theorem}
	To prove this theorem, we recall some useful results. 
\section{Auxiliary results}\label{ppsec2}
	
\subsection{Some definitions and properties}
This section is dedicated to the definition of the concepts.
\begin{definition}[Generalized Lucas sequence]$\\$
The generalized Lucas sequence $\left( \mathcal{U}_n \right)_{n \geqslant 0}$ is defined with initial values $\mathcal{U}_0=0, \ \mathcal{U}_1=1$ and the libear recurrence 
	$$ \mathcal{U}_{n+1}=r \mathcal{U}_n +s \mathcal{U}_{n-1}$$
	where $r$ and $s$ are integers such that $\Delta = r^2 +4 s >0$.
\end{definition}
\noindent For more information about this sequence, the reader can refer to the book of Ribenboim (My Numbers, my friends).
\begin{definition}[Binet formulat]$\\$
The Binet's formula of the generalized Lucas sequence is given by
	$$\mathcal{U}_n=\cfrac{\delta^n - \gamma^n
	}{\delta - \gamma}$$
	where $\delta=\cfrac{r+\sqrt{\Delta}}{2}$ and $\gamma=\cfrac{r-\sqrt{\Delta}}{2}$.
\end{definition}

\begin{definition}[Repdigits]$\\$
A repdigit or sometimes monodigit is a natural number composed of repeated instances of the same digit in a positional number system (often implicitly decimal).
\end{definition}
\noindent The word is a portmanteau of repeated and digit. Examples are $11, 666, 4444$, and $999999$. All repdigits are palindromic numbers and are multiples of repunits. Other well-known repdigits include the repunit primes and in particular the Mersenne primes (which are repdigits when represented in binary).

\begin{definition}[Repdigits are the representation in base $b$]$\\$
Repdigits are the representation in base $b$ of the number $a\left( \dfrac{b^k-1}{b-1}\right)$ where $1\leqslant a\leqslant b-1$  is the repeated digit and $1< k$ is the number of repetitions. 
\end{definition}
\noindent For example, the repdigit $999999$ in base $10$ is $\displaystyle 9\times \dfrac{10^6-1}{10-1}$.

\noindent Recently, the following result was proved in \cite{adedji3}.
\begin{lemma}\label{Lem1}$\\$
		The $n$-th term of the generalized Lucas sequence $\left(\mathcal{U}_n\right)_{n \geqslant 0}$, with $s \in \{-1,1\}$, satisfies the inequalities 
		$$ \delta^{n-2} \leqslant \mathcal{U}_n \leqslant \delta^n $$
		for $n \geqslant 2.$
\end{lemma}
\subsection{A lower bound for linear forms in logarithms}
	
The next tools are related to the transcendental approach to solving Diophantine equations. Let $\eta$ be an algebraic number of degree $d,$ let $a_0 >0$ be the leading coefficient of its minimal polynomial over $\bZ$ and let $\eta=\eta^{(1)},\ldots,\eta^{(d)}$ denote its conjugates. The quantity defined by  
\[
 h(\eta)= \frac{1}{d}\left(\log |a_0|+\sum_{j=1}^{d}\log\max\left(1,\left|\eta^{(j)} \right| \right) \right)
\]
is called the logarithmic height of $\eta.$ Some properties of height are as follows. For $\eta_1, \eta_2$ algebraic numbers and $m\in \bZ,$ we have
\begin{align*}
 h(\eta_1 \pm \eta_2) &\leq h(\eta_1)+ h(\eta_2) +\log2,\\
h(\eta_1\eta_2^{\pm 1}) &\leq h(\eta_1) + h(\eta_2),\\
h(\eta_1^m)&=|m|h(\eta_1).
\end{align*}
If $\eta=\dfrac{p}{q}\in\bQ$ is a rational number in reduced form with $q>0,$ then the above definition reduces to $h(\eta)=\log(\max\{|p|,q\}).$
We can now present the famous Matveev result used in this study. Thus, let $\bL$ be a real number field of degree $d_{\bL}$, $\eta_1,\ldots,\eta_s \in \bL$ and $b_1,\ldots,b_s \in \bZ  \setminus\{0\}.$ Let $B\ge \max\{|b_1|,\ldots,|b_s|\}$ and
\[
\Lambda=\eta_1^{b_1}\cdots\eta_s^{b_s}-1.
\]
Let $A_1,\ldots,A_s$ be real numbers with 
\[
A_i\ge \max\{d_{\bL} h(\eta_i), |\log\eta_i|, 0.16\},\quad i=1,2,\ldots,s.
\]
With the above notations, Matveev proved the following result.
\begin{lemma}[Matveev]\label{pplem:Matveev}$\\$
Assume that $\Lambda\neq 0.$ Then
\[
\log|\Lambda|>-1.4\cdot30^{s+3}\cdot s^{4.5}\cdot d_{\bL}^2 \cdot(1+\log d_{\bL})\cdot(1 +\log B)\cdot A_1\cdots A_s.
\]
\end{lemma}

	\subsection{A generalized result of Baker-Davenport}
	\begin{lemma}[Baker-Davenport]\label{pple:Baker-Davenport}
		Assume that $\tau$ and $\mu$ are real numbers and $M$ is a positive integer. Let $p/q$ be the convergent of the continued fraction of the irrational $\tau$ such that $q>6M$, and let $A, \ B,\mu$ be some real numbers with $A>0$ and $B>1$. Let $\varepsilon=||\mu q||-M\cdot||\tau q||$, where $||\cdot||$ denotes the distance from the nearest integer. If    $\varepsilon>0$, then there is no solution to the inequality
		$$
		0 < m \tau - n + \mu <A B^{-k}
		$$
		in positive integers $m$, $n$ and $k$ with 
		$$
		m\leq M\quad \text{ and }\quad k\geq\dfrac{\log(Aq/\varepsilon)}{\log B}.
		$$ 
	\end{lemma}
	\begin{lemma}[G\'uzman, S\'anchez, Luca] \label{GSL}$\\$
		Let $r \geqslant 1$ and $H > 0$ be such that $H > (4r^2)^r$ and
		$H > L/(\log L)^r$. Then
		$$L < 2^r H(\log H)^r.$$
	\end{lemma}
	\section{Proof of the main result}
	In this study, we consider the following Diophantine equation 
	\begin{equation}\label{Eq1}
		a\left(\cfrac{b^k-1}{b-1}\right)=\mathcal{U}_n-\mathcal{U}_m
	\end{equation}
	with $n > m$ and $1<\leqslant a\leqslant b-1.$\\
	From \eqref{Eq1}, we deduce that: $b^{k-1}< \mathcal{U}_n \leqslant \delta^n$, where we used Lemma \ref{Lem1}. So, we get $(k-1) \log b < n \log \delta$ which leads to 
	\begin{equation}\label{R1}
		k< 1+ n \cfrac{\log \delta}{\log b}.
	\end{equation}
	
	Using now Binet's formula for $\left(\mathcal{U}_n \right)_{n \geqslant 0}$, Diophantine equation \eqref{Eq1} becomes:
	$$\cfrac{\delta^n - \gamma^n
	}{\delta - \gamma}-\cfrac{\delta^m - \gamma^m
	}{\delta - \gamma}=a\left(\cfrac{b^k-1}{b-1}\right),$$ 
	which implies that 
	$$\cfrac{\delta^n}{\delta- \gamma} - \cfrac{ab^k}{b-1}= \cfrac{\gamma^m}{\delta- \gamma}+ \cfrac{\delta^n}{\delta- \gamma}-\cfrac{\gamma^m}{\delta- \gamma}-\cfrac{a}{b-1}.$$
	So taking absolute values on both sides we get:
	\begin{equation}\label{Eq2}
		\left| \cfrac{\delta^n}{\delta- \gamma} - \cfrac{ab^k}{b-1} \right|	\leqslant \cfrac{\left| \gamma \right|^m}{\sqrt{\Delta}} + \cfrac{\delta^m}{\sqrt{\Delta}}+\cfrac{\left| \gamma \right|^n}{\sqrt{\Delta}}+\cfrac{a}{b-1}.
	\end{equation}
	Note that $\left| \gamma \right|=\delta^{-1}$. Thus \eqref{Eq2} becomes :
	\begin{align*}
		\left| \cfrac{\delta^n}{\delta- \gamma} - \cfrac{ab^k}{b-1} \right|	 & \leqslant \cfrac{1 }{ \delta^m \sqrt{\Delta}} + \cfrac{\delta^m}{\sqrt{\Delta}}+\cfrac{1 }{ \delta^n \sqrt{\Delta}}+\cfrac{a}{b-1}\\
		& < 3 + \cfrac{\delta^m}{\sqrt{\Delta}}=\cfrac{3 \sqrt{\Delta}+ \delta^m}{\sqrt{\Delta}}.
	\end{align*}
	Since $s \in \{-1,1\}, \quad \delta \geqslant \cfrac{1+\sqrt{5}}{2}$, then 
	\begin{equation}\label{Eq3}
		\left| \cfrac{\delta^n}{\delta- \gamma} - \cfrac{ab^k}{b-1} \right|	 < 3 + \cfrac{\delta^m}{\sqrt{\Delta}} <\cfrac{1+3\sqrt{\Delta}}{\sqrt{\Delta}} \delta^m.	
	\end{equation}
	Dividing both sides of \eqref{Eq3} by $\cfrac{\delta^n}{\sqrt{\Delta}}$, we get
	\begin{align*}
		\left|1- \delta^{-n} b^k\cfrac{a \sqrt{\Delta}}{b-1} \right|	& <  \cfrac{1+3\sqrt{\Delta}}{\sqrt{\Delta}} \cdot \cfrac{\sqrt{\Delta}}{\delta^n}\cdot \delta^m.\\
		& =\cfrac{1+3\sqrt{\Delta}}{\delta^{n-m}}
	\end{align*}
	So we have:
	\begin{equation}\label{Eq3'}
		\left|\Gamma\right|:=\left|1- \delta^{-n} b^k\cfrac{a \sqrt{\Delta}}{b-1} \right| < \cfrac{1+3\sqrt{\Delta}}{\delta^{n-m}}.
	\end{equation}
	Next, we have to show that $\Gamma \neq 0$.\\
	If $\Gamma =0$, we get :
	$$\delta^n=b^k \cfrac{a \sqrt{\Delta}}{b-1}$$
	which leads to 
	$$\delta^{2n}=b^{2k} \cfrac{a^2 \Delta}{(b-1)^2}=x+y\sqrt{\Delta}$$
	for some integers $x$ and $y$ which is a contradiction because $n \geqslant 1$. Thus $\Gamma \neq 0$ and we can apply Matveev result to $\Gamma$.\\
	Now we put
	\[\eta_1=\delta , \quad \eta_2=b, \quad \eta_3=\cfrac{a \sqrt{\Delta}}{b-1}, \]
	\[b_1=-n, \quad b_2=k, \quad b_3=1, \ \text{and}, \ s=3\]
	$L:=\mathbb{Q}(\eta_1,\eta_2,\eta_3)=\mathbb{Q}(\sqrt{\Delta})$ then 
	$$d_L=\left[ \mathbb{Q}(\eta_1,\eta_2,\eta_3: \mathbb{Q}\right]=2.$$
	For the logarithm heights of $\eta_1,\eta_2$ and $\eta_3$, we have :
	$$h \left(\eta_1\right)=\cfrac{1}{2} \log  \delta, \ h \left(\eta_2\right)=\log b $$ and  \begin{align*}
		h \left(\eta_3\right)&= h \left( \cfrac{a \sqrt{\Delta}}{b-1}\right)\leqslant h \left( \cfrac{a}{b-1} + h \left( \sqrt{\Delta}\right) 
		\right)\\
		&\leqslant \log (b-1) + \cfrac{1}{2} \log \Delta \\
		& < \log b + \cfrac{1}{2} \log \Delta.\\
	\end{align*}
	Thus, we can take $A_1=\log \delta , \quad A_2=2\log b$ and $A_3=2 \log b +\log \Delta.$\\
	Applying Matveev's theorem, we have:
	\begin{equation}\label{Eq4}
		\log \left| \Gamma \right| > -1.4 \cdot 30^6 \cdot 3^{4.5} \cdot 2^2 \cdot (1+ \log2)\cdot (1+\log D)  \cdot \log \delta \cdot 2 \log b \cdot (2 \log b + \log \Delta)
	\end{equation}
	where $D=\max \{\left|b_1 \right|, \left|b_2 \right|,\left|b_3 \right|\}=\{ 1,n,k\}.$
	Since $k<1+n \cfrac{\log \delta}{\log b}$ for $b \geqslant 2$, we can take 
	$$D=1+n \cfrac{\log \delta}{\log b} $$
	Combining \eqref{Eq3'} and \eqref{Eq4} , we get:
	\begin{align*}
		(n-m)\log \delta - \log (1+3\sqrt{\Delta}) &< 1.4 \cdot 30^6 \cdot 3^{4.5} 2^2  (1+ \log2) (1+\log D)   \log \delta \cdot 2 \log b \cdot (2 \log b + \log \Delta)\\
		& < 2 \cdot 10^{12} (1+\log D)   \log \delta 2 \log b  (2 \log b + \log \Delta)
	\end{align*}
	We rewrite Diophantine equation \eqref{Eq1} to obtain that 
	\[\cfrac{\delta^n}{\sqrt{\Delta}} - \cfrac{\delta^m}{\sqrt{\Delta}} - \cfrac{a b^k}{b-1}=\cfrac{\gamma^n}{\sqrt{\Delta}} -\cfrac{\gamma^m}{\sqrt{\Delta}} - \cfrac{a}{b-1}.\]
	Taking absolute values on both sides, we have
	\[\left| \cfrac{\delta^n}{\sqrt{\Delta}} - \cfrac{\delta^m}{\sqrt{\Delta}} - \cfrac{a b^k}{b-1} \right| \leqslant \cfrac{1}{\delta^n \sqrt{\Delta}}+\cfrac{1}{\delta^m \sqrt{\Delta}}+ \cfrac{a}{b-1}<3.\]
	So, we have :
	\[\left|\cfrac{\delta^n}{\sqrt{\Delta}} \left(1 - \delta^{m-n}\right) - \cfrac{a b^k}{b-1}\right|<3.\]
	Dividing both sides by $\cfrac{\delta^n}{\sqrt{\Delta}} \left(1 - \delta^{m-n}\right)$, we get that 
	\begin{equation}\label{Eq5}
		\begin{split}
			\left| 1 - \delta^{-n} \cdot b^k \cfrac{a \sqrt{\Delta}}{(b-1)(1-\delta^{m-n})}\right| & < \cfrac{3 \sqrt{\Delta}}{\delta^n(1-\delta^{m-n})}\\
			& = \cfrac{3 \sqrt{\Delta}\cdot \delta^{n-m}}{\delta^n(\delta^{n-m}-1)}.
		\end{split}
	\end{equation}
	Moreover, $n-m \geqslant 1.$ Let us show it.\\
	From equation \eqref{Eq1}, we have $\mathcal{U}_n -\mathcal{U}_m > 0.$ So 
	\[\delta^{m-2} \leqslant \mathcal{U}_m < \mathcal{U}_n < \delta^n.\]
	Then $m-2 < n$ which implies that $n-m \geqslant -1$.\\
	Note that $n-m$ cannot be equal to $-1$ or $0$. Therefore we have to consider \[n-m \geqslant 1.\]
	Since $n-m \geqslant 1$, then $\delta^{n-m}\geqslant \delta \geqslant \alpha= \cfrac{1+\sqrt{5}}{2}.$\\
	Using now the fact that the numerical function  $f(x)=\cfrac{x}{x-1}$ is decreasing for $x\geqslant \cfrac{1+\sqrt{5}}{2}$, we have
	\[ \cfrac{\delta^{n-m}}{\delta^{n-m}-1} \leqslant \cfrac{\alpha}{\alpha-1} < 2.7.\]
	Hence \eqref{Eq5} becomes 
	\begin{equation}\label{Eq6}
		\left| 1 - \delta^{-n} \cdot b^k \cdot \cfrac{a \sqrt{\Delta}}{(b-1)(1-\delta^{m-n})}\right| < \cfrac{8.1 \cdot \sqrt{\Delta}}{\delta^n}.
	\end{equation}
	Now set $\left| \Gamma^{'}\right|:=\left|  1 - \delta^{-n} \cdot b^k  \cfrac{a \sqrt{\Delta}}{(b-1)(1-\delta^{m-n})} \right|$. 

	Similarly, we can show that $\Gamma^{'} \neq 0.$ \\
	Put 
	\[\eta_1=\delta ,\quad, \eta_2=b, \quad \eta_3=\cfrac{a \sqrt{\Delta}}{(b-1)(1-\delta^{m-n})}\]
	\[b_1=-n, \quad b_2= k , \quad b_3=1.\]
	Note that 
	\begin{align*}
		h(\eta_3) & = h(\cfrac{a \sqrt{\Delta}}{(b-1)(1-\delta^{m-n})}) \\
		& \leqslant h(\cfrac{a}{b-1})+ h (\sqrt{\Delta})+ h (\cfrac{1}{1-\delta^{m-n }})\\
		& < \log b + \cfrac{1}{2} \log \Delta
		+ (n-m)\cdot \cfrac{\log 
			\delta}{2}+ \log 2\\
		& = \log(2 b \sqrt{\Delta}) + \cfrac{n-m}{2} \log \delta.
	\end{align*}
	\[h(\eta_3) < \log(2 b \sqrt{\Delta}) + 10^{12} (1+\log D)  \cdot \log \delta \cdot  \log b \cdot (2 \log b + \log \Delta) + \cfrac{\log (1+3 \sqrt{\Delta})}{2}.\]
	Thus, we can take 
	\begin{align*}
		A_3 & =2 \log(2 b \sqrt{\Delta}) +2\cdot 10^{12} (1+\log D)  \cdot \log \delta \cdot  \log b \cdot (2 \log b + \log \Delta) + \log (1+3 \sqrt{\Delta})\\
		&=  \log(4 b^2 \Delta (1+3 \sqrt{\Delta}) +2\cdot 10^{12} (1+\log D)  \cdot \log \delta \cdot  \log b \cdot (2 \log b + \log \Delta).
	\end{align*}
	By Matveev, we get that :
	\[\log \left|\Gamma^{'}\right| > - 1.4 \cdot 30^6 \cdot 3^{4.5}(1 + \log 2)(1+\log D)\cdot \log \delta \cdot(2 \log b) \cdot A_3. \]
	Combining this with \eqref{Eq6}, we have
	\begin{equation}\label{R2}
		n \log \delta - \log (8.1 \cdot \sqrt{\Delta})< 2 \cdot 10^{12} (1+\log D)\cdot \log \delta \cdot \log b \cdot A_3.	
	\end{equation}
	From \eqref{R1} and \eqref{R2}, we have the proof of Theorem \ref{main}.
	\section{Application: Pell numbers in decimal base}
	In this section, we explicitly determine all repdigits which can be written as difference of two Pell numbers. So our result in this case is Theorem \ref{theo}. 
	In this case, $\mathcal{U}_n$ is Pell number. We have :
	\[(r,s)=(2,1), \quad \Delta=8, \quad \text{and} \  \delta=1+\sqrt{2}.\]
	By the main theorem \ref{main}, we have:
	\[ n \log (1+\sqrt{2}) - \log (8.1 \cdot \sqrt{8})< 2 \cdot 10^{12} (1+\log 8)\cdot \log (1+\sqrt{2}) \cdot \log 10 \cdot \xi\]
	with \[\xi= \log(4 \times 10^2 \cdot 8 (1+3 \sqrt{8}) +2\cdot 10^{12} (1+\log D)  \cdot \log (1+\sqrt{2}) \cdot  \log 10\cdot (2 \log 10 + \log 8) \]
	and
	\begin{align*}
		D=1+ n\dfrac{\log(1+\sqrt{2})}{\log 10} & < 1+1.3 n\\
		& < 2n \quad \text{for} \ n \geqslant2.
	\end{align*} 
	First,
	\begin{align*}
		\xi & < 10.4 + 2.8 \cdot 10^{13} (1+ \log 2n)\\
		& < 3 \cdot 10^{13} (1+ \log 2n) \quad \text{for} \ n \geqslant 2.
	\end{align*}
	So, we get : 
	\begin{align*}
		n &< 1.4 \cdot 10^{26} (1+\log 2n)^2\\
		&=1.4 \cdot 10^{26} (1+\log 2 + \log n)^2. 
	\end{align*}
	Since $n \geqslant 2,$ we obtain
	\[n < 1.7 \cdot 10^{27} . \log^2 n.\]
	Now, we can apply the Lemma of G\'uzman, S\'anchez, and Luca (Lemma \ref{GSL}) by putting
	\[l=2, \quad L=n, \quad \text{and} \ H=1.7 \cdot 10^{27}.\]
	So, we have: $n < 2^2 \cdot 1.7 \cdot 10^{27} \cdot\left(\log(1.7 \cdot 10^{27})\right)^2$, so
	\[n< 2.7 \cdot 10^{31}.\]
	Next, we need to reduce the bound on $n$ by using the Baker-Davenport reduction method due to Dujella and Peth\"o.\\
	Let 
	\[\Lambda_1=-n \log \delta + k \log 10 + \log (\cfrac{a\sqrt{8}}{9}).\]
	The inequality \eqref{Eq3'} can be written as 
	\[\left| \mathrm{e}^{\Lambda_1}-1  \right|< \cfrac{1+3\sqrt{8}}{\delta^{n-m}}.\]
	Observe that $\Lambda_1 \neq 0$ as $e^{\Lambda_1}-1=\Gamma\neq0.$\\
	Assume that $n-m \geqslant 4 $, then \[\left| \mathrm{e}^{\Lambda_1}-1  \right|< \cfrac{1+3\sqrt{8}}{\delta^{n-m}}<\cfrac{1}{2}.\]
	This implies that:
	\[\left| \Lambda_1  \right|<2 \cfrac{1+3\sqrt{8}}{\delta^{n-m}}\]
	since $ \left| x \right|< 2 \left| \mathrm{e}^{x}-1  \right| $ for every real $x$ with $\left|x\right|<\cfrac{1}{2}.$\\
	Dividing both sides by $\log \delta $, we get that:
	\[\left| k \cfrac{\log 10}{\log \delta}-n + \cfrac{\log (a \sqrt{8}/9)}{\log \delta} \right|< \cfrac{21.6}{\delta^{n-m}}.\] 
	Thus, we can take:
	\[\tau= \cfrac{\log 10}{\log \delta}, \quad \mu =\cfrac{\log (a \sqrt{8}/9)}{\log \delta}, \quad A=21.6, \quad B=\delta=1+\sqrt{2} \quad \omega=n-m. \]
	Moreover $k < 2 n < 5.4 \cdot 10^{31}.$ Then we take $M:=5.4 \cdot 10^{31}.$
	With \textit{Mathematica}, we have 
	$q_{73}=1 189 285 833 530 929 228 438 091 844 076 539, \quad \epsilon=0.0049271 , \quad \text{and} \ n-m \leqslant 95.$\\
	Put now
	\[\Lambda_2=-n \log \delta +k \log 10 + \log \left(\cfrac{a \sqrt{8}}{9\left(1-\delta^{m-n}\right)}\right)\]
	So, the inequality \eqref{Eq6} can be written as 
	\[\left| \mathrm{e}^{\Lambda_2}-1  \right|< \cfrac{8.1\sqrt{8}}{\delta^{n}}.\]
	Note also that $\Lambda_2 \neq 0$ as  $\mathrm{e}^{\Lambda_2}-1=\Gamma^{'}\neq0.$\\
	Assuming  $n \geqslant 5 $, we get 
	\[\left| \mathrm{e}^{\Lambda_2}-1  \right|< \cfrac{8.1\sqrt{8}}{\delta^{n}}<\cfrac{1}{2}\]
	and then 
	\[\left| \Lambda_2 \right|<\cfrac{16.2\sqrt{8}}{\delta^{n}}\]
	Dividing both sides by $\log \delta $, we get that:
	\[\left| k \cfrac{\log 10}{\log \delta}-n + \cfrac{\log (a \sqrt{8}/(9(1-\delta^{m-n})))}{\log \delta} \right| < \cfrac{52}{\delta^n}.\]
	To apply Dujella and Peth\"o result, we can set \[\tau= \cfrac{\log 10}{\log \delta}, \quad \mu =\log \left(\cfrac{a \sqrt{8}}{9(1-\delta^{m-n})}\right), \quad A=52, \quad B=\delta=1+\sqrt{2} \quad \omega=n. \]
		Since $k < 1 +n \cfrac{\log \delta}{\log b}< 2 n < 5.4 \cdot 10^{31}.$ Then we take $M:=5.4 \cdot 10^{31}.$ With \textit{Mathematica}, we get
	$q_{73}=1 189 285 833 530 929 228 438 091 844 076 539, \quad \epsilon=0.429295, \quad \text{and}  \ n \leqslant 91.$ Therefore we have proved the Theorem \ref{theo}.
	
	\section*{Acknowledgements}
	The first author is partially supported by Universit\'e de Kara (Togo). The authors thank  K. N Ad\'edji of the Institut de Math\'ematiques et de Sciences Physiques de l'Universit\'e d'Abomey-Calavi (IMSP) for his remarks and suggestions which considerably improved the quality of this paper. 
	
	\section*{Declarations}
	\textbf{Conflict of interest} There is no conflict of interest related to this paper or this submission. The authors have freely chosen this journal for publication without any consideration.
	
	\section*{Data availability} 
	Not applicable.
	\section*{Author contributions} 
	The first draft of the manuscript was written by Pagdame Tiebekabe and all authors commented on previous versions of the manuscript. All authors read and approved the final manuscript.


%
\end{document}